\documentclass[letterpaper,12pt]{article}
\usepackage{amsmath,amssymb,amsfonts,epsf,epsfig}
\newtheorem{theo}{Theorem}[section]

\newtheorem{lemma}[theo]{Lemma}

\newtheorem{conj}[theo]{Conjecture}

\newenvironment{proof}{\noindent {\sc Proof}.}
                {\phantom{a} \hfill \framebox[2.2mm]{ } \bigskip}

\newenvironment{proofof}{\noindent {\sc Proof of Theorem}}
                {\phantom{a} \hfill \framebox[2.2mm]{ } \bigskip}

\newcommand{\NN}{\mathbb{N}}
\newcommand{\ZZ}{\mathbb{Z}}

\def\int{{\rm int}}

\renewcommand{\phi}{\varphi}

\title{On the Spouse-Loving Variant \\ of the Oberwolfach Problem}
\author{Noah Bolohan and Iona Buchanan \\
{\small University of Ottawa} \\ \\
Andrea Burgess \\ {\small University of New Brunswick} \\ \\
Mateja \v{S}ajna\footnote{Email: msajna@uottawa.ca.  Mailing address: Department of Mathematics and Statistics, University of Ottawa,  Ottawa, ON, K1N 6N5, Canada.}\\
{\small University of Ottawa}
\\ \\
Ryan Van Snick \\ {\small University of New Brunswick} \\}

\begin{document}
\maketitle \baselineskip 17pt

\begin{abstract}
We prove that $K_n+I$, the complete graph of even order with a 1-factor duplicated, admits a decomposition into 2-factors, each a disjoint union of cycles of length $m \ge 5$ if and only if $m|n$, except possibly when $m$ is odd and $n=4m$. In addition, we
show that $K_n+I$ admits a decomposition into 2-factors, each a disjoint union of cycles of lengths $m_1,\ldots,m_t$, whenever $m_1,\ldots,m_t$ are all even.

\medskip
\noindent {\em Keywords:} Oberwolfach Problem, 2-factorization, complete graph plus a 1-factor, resolvable minimum covering by cycles, Spouse-Loving Variant.
\end{abstract}

\section{Introduction}

The well-known {\em Oberwolfach Problem} asks the following: Given  $t$ round tables of sizes $m_1,\ldots,m_t$ such that $m_1+\ldots+m_t=n$, is it possible to seat $n$ people around the $t$ tables for an appropriate number of meals so that every person sits next to every other person exactly once? In graph-theoretic terms, the question is asking whether  $K_n$ can be decomposed into 2-factors, each  a disjoint union of cycles of lengths $m_1,\ldots,m_t$, whenever $m_1+\ldots+m_t=n$. Since a graph with odd-degree vertices cannot admit a 2-factorization, Huang, Kotzig, and Rosa \cite{HuaKot}  proposed the analogous problem for $K_n-I$, the complete graph of even order with a 1-factor removed. They  called it the {\em Spouse-Avoiding Variant} since it models a sitting arrangement of $\frac{n}{2}$ couples, where each person gets to sit next to every other person, except their spouse, exactly once, and never gets to sit next to their spouse.

For both of these two basic variants of the Oberwolfach Problem,  the cases with uniform cycle length were completely solved decades ago \cite{AlsHag, AlsSch, HofSch}. In addition, many solutions are now known for  variable cycle lengths; most notably, the problem is solved for $m_1,\ldots,m_t$ all even \cite{BryDan}; for $t=2$ \cite{Tra}; and for $n \le 40$ \cite{DezFra}. In general, however, it is still open.

The Spouse-Avoiding Variant of the Oberwolfach Problem can also be viewed as the {\em maximum packing}  variant. This paper, however,  pertains to the {\em minimum covering} version of the problem; in other words, we are interested in decomposing $K_n+I$ (the complete graph of even order with a 1-factor duplicated) into 2-factors, each  a disjoint union of cycles of lengths $m_1,\ldots,m_t$, where $m_1+\ldots+m_t=n$. We denote this problem by OP$^+(m_1,\ldots,m_t)$, or OP$^+(n;m)$ when $m_1=\ldots=m_t=m$.
This variant, nicknamed the {\em Spouse-Loving Variant}, models a situation where we wish for each person to sit next to exactly one other person --- their spouse --- twice, and next to every other person exactly once.

Necessary and sufficient conditions for OP$^+(n;3)$ to have a solution were previously determined by Assaf, Mendelsohn, and Stinson \cite{AssMenSti}, and Lamken and Mills \cite{LamMil}
under the term {\em resolvable minimum coverings by triples}.

\begin{theo}\cite{AssMenSti,LamMil}\label{thm:AssMenSti}
OP$^+(3t;3)$ has a solution if and only if $t$ is even and $t \ge 6$.
\end{theo}

In this paper, we prove the following main result.

\begin{theo}\label{the:odd}
Let $n$ be an even integer and $m \ge 5$ an integer such that $n \ne 4m$ when $m$ is odd. Then OP$^+(n;m)$ has a solution if and only if $m|n$.
\end{theo}

In addition, we show that the complete solution to the Oberwolfach Problem for bipartite 2-factors \cite{AlsHag,BryDan} implies the following.

\begin{theo}\label{the:even}
Let $m_1,\ldots,m_t$ be even integers greater than 3. Then OP$^+(m_1,\ldots,m_t)$ has a solution.
\end{theo}

\section{Preliminaries}

As usual, $K_n$ denotes the complete graph with $n$ vertices, and for $n$ even, $K_n-I$ and $K_n+I$ denote the complete graph $K_n$ with the edges of the 1-factor $I$ deleted and duplicated, respectively. By $K_{\alpha[k]}$ we denote the complete equipartite graph with $\alpha$ parts of size $k$. A cycle of length $m$ (or $m$-cycle) is denoted by $C_m$.

A set $\{ H_1,\ldots, H_k \}$ of subgraphs of a graph $G$ is called a {\em decomposition} of $G$ if $\{ E(H_1),\ldots,$ $ E(H_k) \}$ is a partition of $E(G)$. If this is the case, we write $G=H_1 \oplus \cdots \oplus H_k$.

A {\em 2-factor} in a graph $G$ is a spanning 2-regular subgraph of $G$. A 2-factor consisting of disjoint cycles of lengths $m_1, m_2,\ldots, m_t$, respectively, is called a {\em
$(C_{m_1},C_{m_2},\ldots,C_{m_t})$-factor}, and a
$(C_{m},C_{m},\ldots,C_{m})$-factor is also called a $C_m$-factor.

A {\em 2-factorization} of a graph $G$ is a decomposition of $G$ into 2-factors. A {\em $(C_{m_1},C_{m_2},\ldots,C_{m_t})$-factorization} is a  2-factorization consisting of $(C_{m_1},C_{m_2},\ldots,C_{m_t})$-factors; a $C_{m}$-factorization is defined analogously.

Apart from direct constructions, the following two previous results will form the most important tools for the proofs of our two main theorems.

\begin{theo}\cite{Liu}\label{thm:Liu}
Let $\alpha$, $k$, and $\ell$ be positive integers, $\alpha \ge 2$, and $3 \le \ell \le \alpha k$. Then $K_{\alpha[k]}$ admits a $C_{\ell}$-factorization if and only if
\begin{enumerate}
\item $\alpha k \equiv 0 \pmod{\ell}$,
\item $k(\alpha-1)$ is even,
\item $\ell$ is even if $\alpha=2$, and
\item $(k,\alpha,\ell) \not\in \{ (2,3,3),(6,3,3),(2,6,3),(6,2,6) \}$.
\end{enumerate}
\end{theo}

\begin{theo}\cite{AlsHag,BryDan} \label{thm:BryDan}
Let $m_1,\ldots,m_t$ be even integers greater than 3. Then OP$(m_1,\ldots,m_t)$ has a solution, that is, $K_n-I$ admits a $(C_{m_1},C_{m_2},\ldots,C_{m_t})$-factorization for $n=m_1+\ldots +m_t$.
\end{theo}

\section{Results}

Our Theorem~\ref{the:even} is an easy corollary to Theorem~\ref{thm:BryDan} \cite{BryDan}.

\bigskip

\begin{proofof} \ref{the:even}.
Let $n=m_1+\ldots +m_t$, and let $I$ denote a chosen 1-factor of $K_n$. By Theorem~\ref{thm:BryDan} \cite{BryDan}, there exists a $(C_{m_1},C_{m_2},\ldots,C_{m_t})$-factorization of $K_n-I$. It is easy to see that, since all $m_i$ are even, there exists a 1-factor $I'$ of $K_n$ such that $I \oplus I'$ is a $(C_{m_1},C_{m_2},\ldots,C_{m_t})$-factor. Thus $K_n \oplus I'$ admits a $(C_{m_1},C_{m_2},\ldots,C_{m_t})$-factorization, and OP$^+(m_1,\ldots,m_t)$ has a solution.
\end{proofof}

We now turn our attention to odd-length cycles.

%\vspace{-5mm}

\begin{figure}[t!]
\centerline{\includegraphics[scale=0.75]{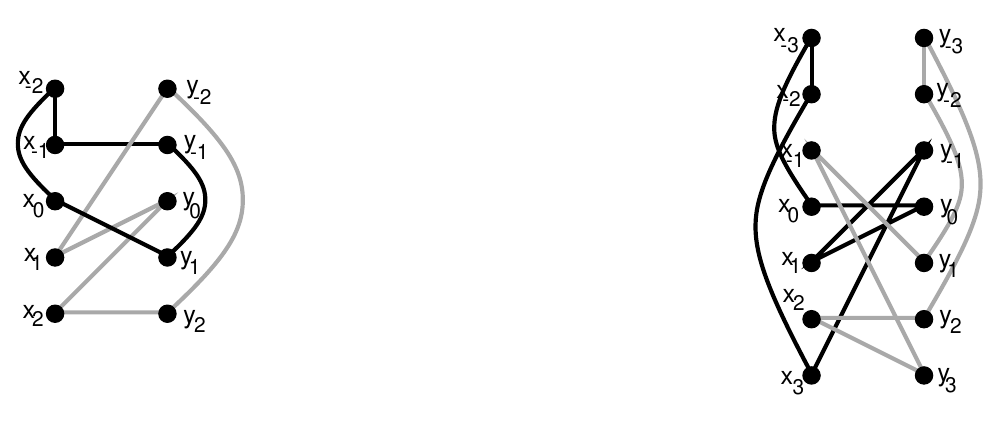}}
\caption{Starter 2-factors in a $C_5$-factorization of $K_{10}+I$ (left) and a $C_7$-factorization of $K_{14}+I$ (right).}\label{fig:m=5and7}
\end{figure}

\begin{lemma}\label{lem:odd}
For all odd $m \ge 5$, there exists a $C_m$-factorization of $K_{2m}+I$.
\end{lemma}

\begin{proof}
Let the vertex set of $K_{2m}+I$ be
$V=\{ x_i: i \in \ZZ_m \} \cup \{ y_i: i \in \ZZ_m \}.$
Edges of the form $x_i x_{i+d}$ and $y_i y_{i+d}$ will be called of {\em pure left} and {\em pure right difference $d$}, respectively, while an edge of the form $x_i y_{i+d}$ will be called of {\em mixed difference $d$}.

Define a permutation $\rho=(x_0 \; x_1\ldots x_{m-1}) (y_0 \; y_1\ldots y_{m-1})$. In each of the cases below, we construct a starter $C_m$-factor $F$ containing exactly one edge of each pure left and pure right difference,  exactly two edges of a single mixed difference, and exactly one edge of every other mixed difference. It then easily follows that $\{ \rho^{\ell}(F): \ell \in \ZZ_m\}$ is a $C_m$-factorization of $K_{2m}+I$.
Figure~\ref{fig:m=5and7} shows the starter 2-factors in a $C_5$-factorization of $K_{10}+I$ and a $C_7$-factorization of $K_{14}+I$, respectively.

\bigskip

\begin{figure}
\centerline{\includegraphics[scale=0.75]{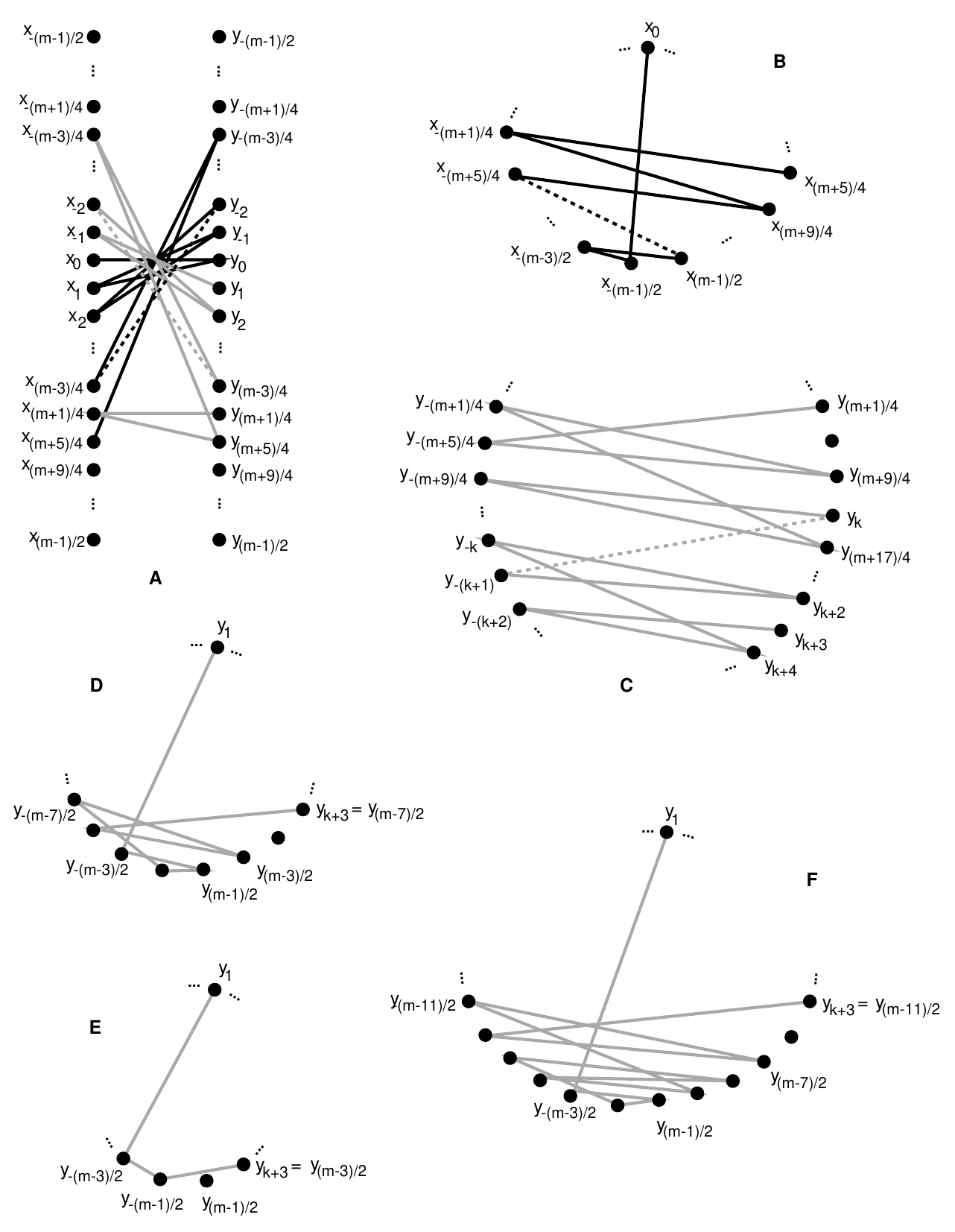}}
\caption{The general construction for Case $m \equiv 3 \pmod{4}$. Paths $P$ and $P'$ (A); path $Q$ (B); the periodic segment of path $Q'$ (C); the non-periodic segment of path $Q'$ for $m \equiv 3,7,11 \pmod{12}$ (D,E,F, respectively). }\label{fig:general3mod4}
\end{figure}

{\sc Case $m \equiv 3 \pmod 4$.}
%Except for $m=11$ (see below), the duplicated 1-factor will be
%$I=\{ x_iy_i: i \in \ZZ_m \}.$
First assume $m \ne 11$. Define the following walks (Figure~\ref{fig:general3mod4}).
\begin{eqnarray*}
P &=& x_0 y_0 x_1 y_{-1} x_2 y_{-2} \ldots x_{\frac{m-3}{4}} y_{-\frac{m-3}{4}} x_{\frac{m+5}{4}}, \\
P' &=& y_1 x_{-1} y_2 x_{-2} \ldots y_{\frac{m-3}{4}} x_{-\frac{m-3}{4}} y_{\frac{m+5}{4}} x_{\frac{m+1}{4}} y_{\frac{m+1}{4}},  \\
Q &=& x_{\frac{m+5}{4}} x_{-\frac{m+1}{4}} x_{\frac{m+9}{4}} x_{-\frac{m+5}{4}} \ldots x_{\frac{m-1}{2}} x_{-\frac{m-3}{2}} x_{-\frac{m-1}{2}} x_0.
\end{eqnarray*}
It is not difficult to verify that each of these walks is in fact a path, and that these three paths are pairwise vertex-disjoint, except that $P$ and $Q$ have the same pair of endpoints.

Observe that $P$ contains exactly one edge of each of the mixed differences in
$$D(P)=\left\{ 0,-1,-2,-3,\ldots,-\frac{m-3}{2}, \frac{m-1}{2} \right\},$$
$P'$ contains exactly one edge of each of the mixed differences in
$$D(P')=\left\{ 0,1, 2,3,\ldots,\frac{m-3}{2}, -\frac{m-1}{2} \right\},$$
while $Q$ contains exactly one edge of each pure left difference in
$$D(Q)=\left\{ 1,2,\ldots, \frac{m-1}{2} \right\}.$$

Our starter 2-factor will contain the $m$-cycle $C=PQ$, as well as the $m$-cycle $C'=P'Q'$, where $Q'$ is the $(y_{\frac{m+1}{4}},y_{1})$-path to be defined as follows below.

First, for $i \le \frac{m-9}{2}$ of the form $i=\frac{m+1}{4}+3j$, where $j \in \NN$, define a 6-path
$$Q_i= y_i y_{-(i+1)} y_{i+2} y_{-i} y_{i+4} y_{-(i+2)} y_{i+3}.$$
Observe that $Q_i$ contains exactly one edge of each of the pure right differences in
$$D(Q_i)=\left\{ -(2i+1), -(2i+2),\ldots,-(2i+6) \right\}.$$
Thus, if $j>0$, paths $Q_i$ and $Q_{i+3j}$ have no edges of the same difference, and share no vertices except that the terminus of $Q_i$  equals the
source of $Q_{i+3}$. We use paths $Q_i$ to define $Q'$ as follows.

If $m \equiv 3 \pmod {12}$, let
$$Q'= Q_{\frac{m+1}{4}} Q_{\frac{m+13}{4}} Q_{\frac{m+25}{4}} \ldots Q_{\frac{m-13}{2}} y_{\frac{m-7}{2}} y_{-\frac{m-5}{2}} y_{\frac{m-3}{2}} y_{-\frac{m-7}{2}} y_{-\frac{m-1}{2}} y_{\frac{m-1}{2}}  y_{-\frac{m-3}{2}} y_1.$$
%Note that for $m=15$, we simply have $Q'=y_{\frac{m-7}{2}} y_{-\frac{m-5}{2}} y_{\frac{m-3}{2}} y_{-\frac{m-7}{2}} y_{-\frac{m-1}{2}} y_{\frac{m-1}{2}}  y_{-\frac{m-3}{2}} y_1.$

If $m \equiv 7 \pmod {12}$, let
$$Q'= Q_{\frac{m+1}{4}} Q_{\frac{m+13}{4}} Q_{\frac{m+25}{4}} \ldots Q_{\frac{m-9}{2}} y_{\frac{m-3}{2}} y_{-\frac{m-1}{2}} y_{-\frac{m-3}{2}}  y_1.$$
%Note that for $m=7$, we simply have $Q'= y_{\frac{m-3}{2}} y_{-\frac{m-1}{2}} y_{-\frac{m-3}{2}}  y_1.$

If $m \equiv 11 \pmod {12}$, let
$$Q'= Q_{\frac{m+1}{4}} Q_{\frac{m+13}{4}}  \ldots Q_{\frac{m-17}{2}} y_{\frac{m-11}{2}} y_{-\frac{m-9}{2}} y_{\frac{m-7}{2}} y_{-\frac{m-11}{2}} y_{\frac{m-3}{2}} y_{-\frac{m-5}{2}} y_{\frac{m-5}{2}} y_{-\frac{m-7}{2}}
y_{-\frac{m-1}{2}} y_{\frac{m-1}{2}}  y_{-\frac{m-3}{2}} y_1.$$
Note that for $m \in \{ 7,15,23\}$, the walk $Q'$ contains no paths $Q_i$.

It is not difficult to verify that in each case, $Q'$ is in fact a path, and contains exactly one edge of each pure right difference in
$$D(Q')=\left\{ 1,2,\ldots,\frac{m-1}{2} \right\}.$$

Moreover, the paths $P'$ and $Q'$ are internally vertex-disjoint and share the endpoints, so $C'=P'Q'$ is an $m$-cycle. Furthermore, the $m$-cycles $C$ and $C'$ are disjoint, and $F=\{ C,C'\}$ is an $C_m$-factor in $K_{2m}+I$ containing exactly one edge of each pure left difference, each pure right difference, and each mixed difference, except that it contains exactly two edges of mixed difference 0. Hence $\{ \rho^{\ell}(F): \ell \in \ZZ_m\}$ is a $C_m$-factorization of $K_{2m}+I$.

For $m=11$, we define $P$ and $Q$ as above, and modify $P'$ and $Q'$ as follows:
\begin{eqnarray*}
P' &=& y_1 x_{-1} y_2 x_{-2} y_4 x_3 y_5, \\
Q' &=& y_5 y_{-4} y_{-5} y_3 y_{-3} y_1.
\end{eqnarray*}
We obtain a $C_m$-factorization of $K_{2m}+I$ (for $m=11$) as before. Note that, in this case, the duplicated mixed difference is $2$.

\bigskip

\begin{figure}
\centerline{\includegraphics[scale=0.75]{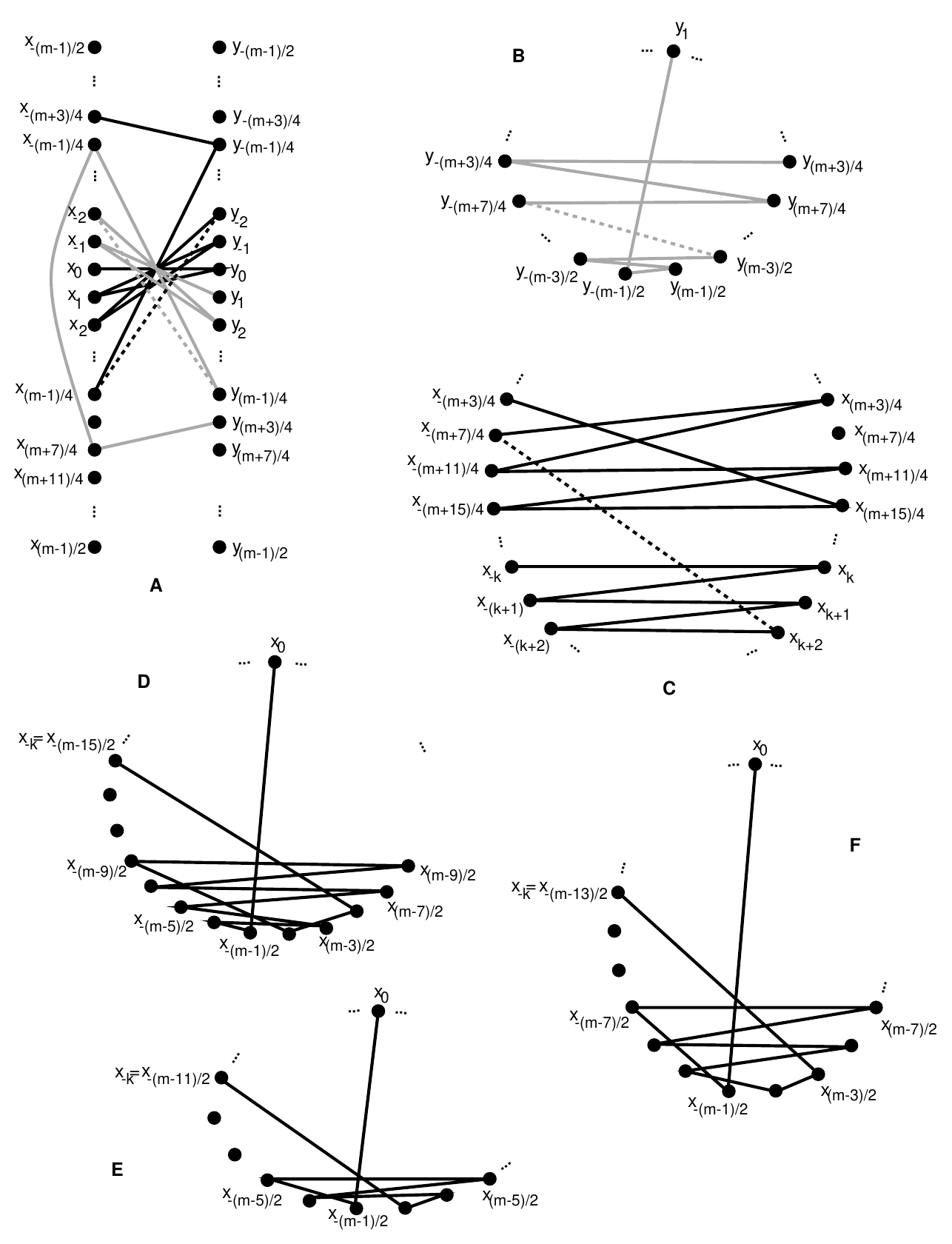}}
\caption{The general construction for Case $m \equiv 1 \pmod{4}$. Paths $P$ and $P'$ (A); path $Q'$ (B); the initial and periodic segments of path $Q$ (C); the non-periodic segment of path $Q$ for $m \equiv 1,5,9 \pmod{12}$ (D,E,F, respectively). }\label{fig:general1mod4}
\end{figure}

{\sc Case $m \equiv 1 \pmod{4}$.} %Except for $m=9$ (see below), the duplicated 1-factor will be $I=\{ x_iy_{-1}: i \in \ZZ_m \}.$
First assume $m \ge 13$. Define the following walks (Figure~\ref{fig:general1mod4}).
\begin{eqnarray*}
P &=& x_0 y_0 x_1 y_{-1} x_2 y_{-2} \ldots x_{\frac{m-1}{4}} y_{-\frac{m-1}{4}} x_{-\frac{m+3}{4}}, \\
P' &=&  y_1 x_{-1} y_2 x_{-2} \ldots  y_{\frac{m-1}{4}} x_{-\frac{m-1}{4}} x_{\frac{m+7}{4}} y_{\frac{m+3}{4}}, \\
Q' &=& y_{\frac{m+3}{4}} y_{-\frac{m+3}{4}} y_{\frac{m+7}{4}} y_{-\frac{m+7}{4}}\ldots y_{\frac{m-3}{2}} y_{-\frac{m-3}{2}}  y_{\frac{m-1}{2}} y_{-\frac{m-1}{2}} y_1.
\end{eqnarray*}

It is not difficult to verify that each of these walks is in fact a path, and that these three paths are pairwise vertex-disjoint, except that $P'$ and $Q'$ have the same endpoints.

Observe that $P$ contains exactly one edge of each of the mixed differences in
$$D(P)=\left\{ 0,-1,-2,\ldots,-\frac{m-1}{2}, 1 \right\},$$
$P'$ contains an edge of pure left difference $\frac{m-3}{2}$ and exactly one edge of each of the mixed differences in
$$D(P')=\left\{ -1,2,3,\ldots,\frac{m-1}{2} \right\},$$
and $Q$ contains exactly one edge of each pure right difference in
$$D(Q')=\left\{ 1,2,\ldots, \frac{m-1}{2} \right\}.$$

Next, let $C'=Q'P'$ and observe that $C'$ is an $m$-cycle. The second cycle in our starter 2-factor will be $C=PQ$, where $Q$ is the $(x_{-\frac{m+3}{4}},x_0)$-path to be defined as follows below.

First, for $i \le \frac{m-11}{2}$ of the form $i=\frac{m+19}{4}+3j$, where $j \in \NN$,  define a 6-path
$$Q_i= x_{-(i-3)} x_{i+2} x_{-(i+2)} x_{i+1} x_{-(i+1)} x_i x_{-i} .$$
Observe that $Q_i$ contains exactly one edge of each pure left difference in
$$D(Q_i)=\{ -(2i-1), -2i,\ldots,-(2i+4) \}.$$
Thus, for $j>0$, paths $Q_i$ and $Q_{i+3j}$  contain no edges of a common difference, and share no vertices except that the terminus of $Q_i$ equals the source of $Q_{i+3}$.

For $m \ge 17$, we additionally  define
$$Q_{\frac{m+7}{4}}= x_{-\frac{m+3}{4}} x_{\frac{m+15}{4}}
x_{-\frac{m+15}{4}} x_{\frac{m+11}{4}}
x_{-\frac{m+11}{4}} x_{\frac{m+3}{4}}
x_{-\frac{m+7}{4}}.$$
Thus $Q_{\frac{m+7}{4}}$ contains exactly one edge of each pure left difference in
$$D(Q_{\frac{m+7}{4}})=\left\{ \frac{m-15}{2},\frac{m-13}{2},\ldots,\frac{m-5}{2}  \right\}$$
and shares no differences and no vertices with any path $Q_{\frac{m+19}{4}+3j}$ for $j \in \NN$, except that the terminus of $Q_{\frac{m+7}{4}}$ equals the source of $Q_{\frac{m+19}{4}}$. We use the paths $Q_i$ to define $Q$ as follows.

If $m \equiv 1 \pmod {12}$ and $m \ge 37$, let
$$Q= Q_{\frac{m+7}{4}} Q_{\frac{m+19}{4}} \ldots Q_{\frac{m-15}{2}}
x_{-\frac{m-15}{2}} x_{\frac{m-5}{2}} x_{\frac{m-1}{2}} x_{-\frac{m-9}{2}} x_{\frac{m-9}{2}} x_{-\frac{m-7}{2}}
x_{\frac{m-7}{2}} x_{-\frac{m-5}{2}}
x_{\frac{m-3}{2}} x_{-\frac{m-3}{2}}
x_{-\frac{m-1}{2}} x_0.$$
Modifications are required for small values of $m$: for $m=13$, we take $Q=x_{-4} x_6 x_4 x_{-5} x_{-6} x_0$, and
for $m=25$, we  take
$Q=x_{-7} x_{10} x_{12} x_{-8} x_{7} x_{-9} x_9 x_{-10} x_{11} x_{-11} x_{-12} x_0$ instead.

If $m \equiv 5 \pmod {12}$ and $m \ge 29$, let
$$Q= Q_{\frac{m+7}{4}} Q_{\frac{m+19}{4}} Q_{\frac{m+31}{4}} \ldots Q_{\frac{m-11}{2}}
x_{-\frac{m-11}{2}} x_{\frac{m-1}{2}} x_{\frac{m-3}{2}} x_{-\frac{m-3}{2}} x_{\frac{m-5}{2}} x_{-\frac{m-5}{2}} x_{-\frac{m-1}{2}} x_0,$$
while for $m=17$, we take
$Q=x_{-5} x_8 x_7 x_{-7} x_{5} x_{-6} x_{-8} x_0.$

If $m \equiv 9 \pmod {12}$ and $m \ge 33$, let
$$Q= Q_{\frac{m+7}{4}} Q_{\frac{m+19}{4}} Q_{\frac{m+31}{4}} \ldots Q_{\frac{m-13}{2}}
x_{-\frac{m-13}{2}} x_{\frac{m-3}{2}} x_{\frac{m-1}{2}} x_{-\frac{m-3}{2}} x_{\frac{m-5}{2}} x_{-\frac{m-5}{2}}
x_{\frac{m-7}{2}} x_{-\frac{m-7}{2}}
x_{-\frac{m-1}{2}} x_0,$$
while for $m=21$, take
$Q=x_{-6} x_9 x_{10} x_{-9} x_{8} x_{-8} x_6 x_{-7} x_{-10} x_0$ instead.

It is not difficult to verify that in each case, the walk $Q$ is in fact a path, and contains exactly one edge of each pure left difference in
$$D(Q)=\left\{ 1,2,3, \ldots,\frac{m-5}{2},\frac{m-1}{2} \right\}.$$
Recall that pure left difference $\frac{m-3}{2}$ was used in the cycle $C'$.

Since paths $P$ and $Q$ are internally vertex-disjoint and share their endpoints, $C=PQ$ is indeed an $m$-cycle. Moreover, $F=\{ C,C'\}$ is a pair of disjoint $m$-cycles in $K_{2m}+I$ containing exactly one edge of each pure left difference, each pure right difference, and each mixed difference, except it contains exactly two edges of mixed difference $-1$. Hence $F$ is a starter 2-factor of a $C_m$-factorization of $K_{2m}+I$.

For $m=9$, we define paths $P$ and $Q'$ as before, and  modify $P'$ and $Q$ as
$$P' = y_1 x_{-1} y_2 x_{-2} x_{-4} y_3 \qquad \mbox{ and } \qquad Q=x_{-3} x_3 x_4 x_0.$$
In this case, the duplicated mixed difference is $-2$.

Finally, for $m=5$, we let
$$C=x_0 y_1 y_{-1} x_{-1} x_{-2}x_0 \qquad \mbox{ and } \qquad  C'=y_0 x_1 y_{-2} y_2 x_2 y_0,$$
so that $F=\{ C,C'\}$ is a starter 2-factor with a duplicated mixed difference $0$.

This completes all cases of the construction of a $C_m$-factorization of $K_{2m}+I$.
\end{proof}

\bigskip

We round up this section with the proof of our main result.

\bigskip

\begin{proofof} \ref{the:odd}.
Let $n$ be an even integer and $m \ge 5$. If OP$^+(n;m)$ has a solution, then clearly $m|n$.

Conversely, assume $n$ is even and $n=tm$ for an integer $t$. If $m$ is even, then OP$^+(n;m)$ has a solution by Theorem~\ref{the:even}. Hence assume $m$ is odd and $t \ne 4$. The result follows from Lemma~\ref{lem:odd} if $t=2$. Hence let $t \ge 6$. Decompose $K_{tm}+I$ into $\frac{t}{2}$ disjoint copies of $K_{2m}+I$ and the complete equipartite graph with $\frac{t}{2}$ parts of size $2m$. The first graph admits a $C_m$-factorization by Lemma~\ref{lem:odd}, and the second by Theorem~\ref{thm:Liu}. Hence $K_{n}+I$ admits a $C_m$-factorization, and OP$^+(n;m)$ has a solution.
\end{proofof}

\section{Conclusion}

The main result of this paper is an almost complete solution to the Spouse-Loving Variant of the Oberwolfach Problem with uniform cycle lengths. The only case that remains unsolved is captured by the following conjecture.

\begin{conj}\label{conj:4m}
Let $m$ be a positive odd integer. Then OP$^+(4m;m)$ has a solution if and only if $m \ge 5$.
\end{conj}

It is worthwhile observing that this case is analogous to the only case left unsolved in \cite{AlsSch}, where the Oberwolfach Problem with uniform cycle lengths was almost completely solved. The reason is the same: the reduction to case $n=2m$ does not work since $K_{2m,2m}$ admits no $C_m$-factorization when $m$ is odd. A solution to $OP(4m;m)$, which appeared in \cite{HofSch} two years after \cite{AlsSch}, required a rather different and highly involved construction.

The reader may be wondering whether a solution to OP$^+(4m;m)$ may be constructed similarly to our solution to OP$^+(2m;m)$. Suppose this is indeed the case. Analogously to the proof of  Lemma~\ref{lem:odd}, denote $V(K_{4m}+I)= \{ x_i: i \in \ZZ_{2m} \} \cup  \{ y_i: i \in \ZZ_{2m} \}$ and $\rho=(x_0 \; x_1 \ldots x_{2m-1})(y_0 \; y_1 \ldots y_{2m-1})$, and let $F$ be a starter 2-factor such that $\{ \rho^i(F): i \in \ZZ_{2m} \}$ is a $C_m$-factorization of $K_{4m}+I$. Then $F$ contains exactly one edge of each difference in $K_{4m}+I$, and since the edges of pure left difference $m$ and pure right difference $m$ lie in orbits of length $m$, they must be the ones forming the duplicated 1-factor $I$. Hence the freedom of choosing $I$ that we enjoyed in the case $n=2m$ is lost. However, the fatal flaw is as follows. The graph $K_{4m}+I$ contains edges of exactly $4m$ distinct differences (which is promising), of which exactly $2m+1$ are odd ($\frac{m+1}{2}$ of the pure left differences, $\frac{m+1}{2}$ of the pure right differences, and $m$ of the mixed differences). However, each cycle in $F$ must contain an even number of odd differences --- a contradiction.

Thus, proving Conjecture~\ref{conj:4m} will require a different approach. So far we have been able to verify it for $5 \le m \le 23$, using a combination of general construction, case-specific construction, and, for larger $m$, a computer search. However, despite considerable effort, we have to conclude that a complete proof of Conjecture~\ref{conj:4m} is presently out of our reach.

\newpage

\centerline{\bf Acknowledgement}

\bigskip

A. Burgess and M. \v{S}ajna gratefully acknowledge support by the Natural Sciences and Engineering Research Council of Canada.

\end{document}